\documentclass[12pt]{article}
\usepackage{amsmath,amssymb}

\def\R{{\Bbb R}}
\def\C{{\Bbb C}}

\def\P{{\Bbb P}}

\begin{document}
\begin{center}
{\LARGE
Almost periodic mappings to complex manifolds.}
\end{center}
\vskip 0.1cm

\centerline{\large S.~Favorov, N.~Parfyonova }
\vskip 0.7cm

{\bf Abstract.} H.Bohr in 1930 proved that if a holomorphic
bounded function  on a strip is almost periodic on a straight line
in the strip, then it is almost periodic on the whole strip. We
find some conditions when the result is valid for holomorphic
mappings of  tube domains to various complex manifolds.

\large

\bigskip

A continuous mapping $F$ of a tube $T_K=\{z=x+iy: x\in\R^m,y\in K\subset\R^m\}$ to a metric space
$X$ is {\it almost periodic} if the family $\{F(z+t)\}_{t\in\R^m}$ of shifts along $\R^m$ is a relatively
compact set with respect to the topology of the uniform convergence on $T_K.$

Further, let $X$ be a complex manifold, and $F$ be a holomorphic mapping of a tube
$T_{\Omega}=\{z=x+iy: x\in\R^m,y\in \Omega\},$ with the convex open base $\Omega\subset \R^m,$ to $X.$
We will say that $F$ is {\it almost periodic} if the restriction of $F$ to each tube $T_K,$ with the
compact base $K\subset\Omega,$ is almost periodic.

For $X=\C$ we obtain the well-known class of holomorphic almost
periodic functions; for $X=\C^m$ the corresponding class was being
studied in \cite{FRR}, \cite{Rashkovskii1}, \cite{Rashkovskii2},
\cite{P3}; for $X=\C\P,$ we get the class of meromorphic almost
periodic functions that was being studied in \cite{PF},
\cite{SiB}, \cite{P2}; the class of holomorphic almost periodic
curves, corresponding to the case $X=\C\P^m,$ was being studied in
\cite{P3}.

The following theorem is well known:

{\bf Theorem B} (H. Bohr \cite{bor}) {\it If a holomorphic bounded function
 on a strip is almost periodic on some straight line in this strip, then this
 function is almost periodic on the whole strip.}

This theorem was extended to holomorphic functions on a tube domain in \cite{Udodova};
besides usual uniform metric, various integral metrics were being studied there.

The direct generalization of Theorem B to complex manifold is not true:

{\bf Example.} Let $f(z)$ be the quotient of two periodic functions $\sin\sqrt{2}z$
and $\sin z.$ It is clear that the restriction of $f(z)$ to any straight line
 $\{z=x+iy_0:x\in \R\}, y_0\ne 0,$ is almost periodic on this line. Besides,
$f(z)$ is bounded as a mapping to the compact manifold $\C\P.$ Nevertheless, zeros
and poles of $f(z)$ are not separated on the real axis $z=x\in \R,$ therefore $f(z)$
is not almost periodic on any substrip containing the real axis ( see \cite{PF}).

In order to give the right version of Theorem B, we need the following proposition.

{\bf Proposition.} {\it  If $K$ is compact subset of $\R^m,$ then an
 almost periodic mapping $F$ of $T_K$ to a metric space $X$ is uniformly continuous
 and $F(T_K)$ is a relatively compact subset of $X$.}

{\bf Proof.} Let
$$ \varphi(t)=\sup\limits_{y\in K} d(F(z+t),F(z)),\quad t\in \R^m.
$$
It is easy to prove that $\varphi(t)$ is an almost periodic
function on $\R^m,$ therefore the function $\varphi(t)$ extends
continuously to Bohr's compactification $B$ of the set $\R^m$
(see, for example, \cite{Berg}). Hence for any point $\tau\in B$
and any $\varepsilon>0$ there exists a neighborhood $U\subset B$
of this point such that

$$|\varphi(t)-\varphi(t')|<\varepsilon
$$
for each $t, t'\in U\cap \R^m.$

Using the obvious inequality
$$\varphi(t-t')\le  |\varphi(t)-\varphi(t')|$$
we get for all $z\in T_K, t, t'\in U\cap \R^m$
$$d(F(z+t'),F(z+t))\le  \varepsilon.$$
Therefore the function
$$\overline{F}(\tau,y)=\lim\limits_{{x\to\tau,}\  {x\in\R^m}}
F(x+iy)$$ is well-defined and continuously maps  the compact set
$B\times K$ to $X.$ Since $\overline{F}(x,y)=F(x+iy)$ for $x\in
\R^m, $ we obtain all the statements of our proposition.

Note that a bounded holomorphic function on a tube domain is
uniformly continuous on every subtube with the compact base, but
bounded holomorphic mappings, in general, have no this property.
Therefore the following result is natural. {\bf Theorem.} {\it Let
$F$ be a holomorphic mapping of $T_{\Omega}$ to a complex manifold
$X$ such that for every compact subset $K\subset\Omega$ the
mapping $F$ is uniformly continuous on $T_K$ and $F(T_K)$ is a
relatively compact subset of $X.$ If the restriction of $F(z)$ to
some hyperplane $\R^m+iy'$ is almost periodic, then $F(z)$ is an
almost periodic mapping of $T_\Omega$ to $X$.}

{\bf Corollary.} {\it Let $F$ be a holomorphic mapping form
$T_\Omega$ to a compact complex manifold $X$ such that $F$ is
uniformly continuous on $T_K$ for every compact set $K\subset
\Omega.$ If the restriction of $F(z)$ to some hyperplane
$\R^m+iy'$ is almost periodic, then $F(z)$ is an almost periodic
mapping of $T_\Omega$ to $X$.} {\bf Proof the Theorem.} Take an
arbitrary sequence $\{t_n\}\subset \R^m.$ Since the function
$F(z)$ is uniformly continuous, the family $\{F(z+t_n\}$ is
equicontinuous on each compact set $S\subset T_{\Omega}.$ Further,
it follows from the condition of the Theorem that the union of all
the images of S under mappings of this family is contained in a
compact subset of $X.$ Therefore, passing on to a subsequence if
necessary, we may assume that the sequence $\{F(z+t_n)\}$
converges  to a holomorphic mapping $G(z)$ uniformly on every
compact subset of $T_{\Omega}.$ It easy to see that the mapping
$G(z)$ is bounded and uniformly continuous on every tube $T_{K}$
with the compact base $K\subset\Omega.$ Let us prove that this
convergence is uniform on every $T_K.$ Assume the contrary. Then
we get
\begin{equation}\label{p_eq1}
d(F(z_n+t_n),G(z_n))\ge \varepsilon_0>0
\end{equation}
for some sequence $z_n=x_n+iy_n\in T_{K'},$ where $K'$ is some
compact subset of $\Omega.$ Replacing sequence by a subsequence if
necessary, we may assume that the mappings $G(x_n+z)$ converge to
a holomorphic mapping ${H}(z),$ and the mappings $F(z+x_n+t_n)$
converge to a holomorphic mapping $\tilde{H}(z)$  uniformly on
every compact subsets of $T_{\Omega}.$ We may also assume that
$y_n\to y_0\in K'.$ Using (\ref{p_eq1}) we get
$$|\tilde{H}(iy_0)-H(iy_0))|\ge\varepsilon_0.$$
Since the mapping $F(x+iy')$ of $\R^m$ to $X$ is almost periodic,
we may assume that a subsequence of the mappings $F(x+t_n+iy')$
converges to  $G(x+iy')$ uniformly in $x\in\R^m.$ Therefore the
sequences of mappings $F(x+x_n+t_n+iy')$ and $G(x+x_n+iy')$ have
the same limit, i. e., $\tilde{H}(x+iy')=H(x+iy')$ for all
$x\in\R^m.$ Since $\tilde{H}(z),$  ${H}(z)$ are holomorphic
mappings, we get $\tilde{H}(x+iy)\equiv H(x+iy)$ on $T_{\Omega}.$
This contradiction proves the Theorem.
\renewcommand{\refname}{

\end{document}